\documentclass{amsart}
\usepackage{amsfonts, amsmath, amssymb, epsfig, booktabs, url}
\usepackage[normalem]{ulem}

\makeatletter
\def\labelenumi{\textnormal{(\@alph\c@enumi)}}
\def\theenumi{\@alph \c@enumi}
\def\labelenumii{\textnormal{(\@roman\c@enumii)}}
\def\theenumii{\@roman \c@enumii}
\newcount\charno
\def\alphapart#1{\charno=96
\advance\charno by#1\char\charno}

\makeatother

\expandafter\ifx\csname BibSpec\endcsname\relax\else
\BibSpec{collection.article}{%
    +{}  {\PrintAuthors}                {author}
    +{,} { \textit}                     {title}
    +{.} { }                            {part}
    +{:} { \textit}                     {subtitle}
    +{,} { \PrintContributions}         {contribution}
    +{,} { \PrintConference}            {conference}
    +{}  {\PrintBook}                   {book}
    +{,} { }                            {booktitle}
    +{,} { }                           {series}
    +{,} { \voltext}                    {volume}
    +{,} { }                            {publisher}
    +{,} { }                            {organization}
    +{,} { }                            {address}
    +{,} { \PrintDateB}                 {date}
    +{,} { pp.~}                        {pages}
    +{,} { }                            {status}
    +{,} { \PrintDOI}                   {doi}
    +{,} { available at \eprint}        {eprint}
    +{}  { \parenthesize}               {language}
    +{}  { \PrintTranslation}           {translation}
    +{;} { \PrintReprint}               {reprint}
    +{.} { }                            {note}
    +{.} {}                             {transition}
    +{}  {\SentenceSpace \PrintReviews} {review}
}
\BibSpec{article}{%
    +{}  {\PrintAuthors}                {author}
    +{,} { \textit}                     {title}
    +{.} { }                            {part}
    +{:} { \textit}                     {subtitle}
    +{,} { \PrintContributions}         {contribution}
    +{.} { \PrintPartials}              {partial}
    +{,} { }                            {journal}
    +{}  { \textbf}                     {volume}
    +{}  { \PrintDatePV}                {date}
    +{,} { \issuetext}                  {number}
    +{,} { \eprintpages}                {pages}
    +{,} { }                            {status}
    +{,} { \PrintDOI}                   {doi}
    +{,} { available at \eprint}        {eprint}
    +{}  { \parenthesize}               {language}
    +{}  { \PrintTranslation}           {translation}
    +{;} { \PrintReprint}               {reprint}
    +{.} { }                            {note}
    +{.} {}                             {transition}
    +{}  {\SentenceSpace \PrintReviews} {review}
}
\BibSpec{book}{%
    +{}  {\PrintPrimary}                {transition}
    +{,} { \textit}                     {title}
    +{.} { }                            {part}
    +{:} { \textit}                     {subtitle}
    +{,} { \PrintEdition}               {edition}
    +{}  { \PrintEditorsB}              {editor}
    +{,} { \PrintTranslatorsC}          {translator}
    +{,} { \PrintContributions}         {contribution}
    +{,} { }                            {series}
    +{,} { \voltext}                    {volume}
    +{,} { }                            {publisher}
    +{,} { }                            {organization}
    +{,} { }                            {address}
    +{,} { pp.~}                        {pages}
    +{,} { \PrintDateB}                 {date}
    +{,} { }                            {status}
    +{}  { \parenthesize}               {language}
    +{}  { \PrintTranslation}           {translation}
    +{;} { \PrintReprint}               {reprint}
    +{.} { }                            {note}
    +{.} {}                             {transition}
    +{}  {\SentenceSpace \PrintReviews} {review}
}
\fi
\newcommand{\bt}{\begin{table}[tp]\centering\footnotesize}
\newcommand{\et}{\end{table}}
\newcommand{\btab}[1]{\begin{tabular}{#1}}
\newcommand{\etab}{\end{tabular}}

\newcommand{\bfig}{\begin{figure}[tp]\centering}
\newcommand{\efig}{\end{figure}}

\newcommand{\beq}{\begin{equation}}
\newcommand{\eeq}{\end{equation}}
\newcommand{\beqnn}{\begin{equation*}}
\newcommand{\eeqnn}{\end{equation*}}
\newcommand{\beqa}{\begin{eqnarray}}
\newcommand{\eeqa}{\end{eqnarray}}
\newcommand{\beqann}{\begin{eqnarray*}}
\newcommand{\eeqann}{\end{eqnarray*}}

\newcommand{\C}{{\mathbb{C}}}

\newcommand{\R}{{\mathbb{R}}}

\renewcommand{\Re}{{\operatorname{Re}}}

\newcommand{\bm}[1]{\mbox{\boldmath{$ #1 $}}}
\newcommand{\sbm}[1]{\bm{\scriptsize{\bm{#1}}}}

\allowdisplaybreaks
\begin{document}
\title[Correction to Chebyshev Approximation Coefficients]{Correction to Partial Fraction Decomposition Coefficients for Chebyshev Rational Approximation on the Negative Real Axis}

\author[M. Pusa]{Maria Pusa}
\address{VTT Technical Research Centre of Finland, Nuclear Energy, P.O. Box 1000, FI-02044 VTT, Finland}
\email{maria.pusa@vtt.fi}

\begin{abstract}
  Chebyshev rational approximation can be a viable method to compute the exponential of matrices with eigenvalues in the vicinity of the negative real axis, and it was recently applied successfully to solving nuclear fuel burnup equations. Determining the partial fraction decomposition ({\scriptsize PFD}) coefficients of this approximation can be difficult and they have been provided (for approximation orders 10 and 14) by Gallopoulos and Saad in ``Efficient solution of parabolic equations by Krylov approximation methods'', SIAM J. Sci. Stat. Comput., 13(1992). It was recently discovered that the order 14 coefficients contain errors and result in $10^2$ times poorer accuracy than expected by theory. The purpose of this note is to provide the correct {\scriptsize PFD} coefficients for approximation orders 14 and 16 and to briefly discuss the approximation accuracy resulting from the erroneous coefficients.
\end{abstract}

\date{12 June 2012}
\maketitle


\section{Chebyshev rational approximation}

This note concerns the computation of matrix exponential based on {the} Chebyshev rational approximation method (abbreviated {\sc cram} in~\cite{pusa2010}) on the negative real axis.
In this approach, the matrix exponential $e^{\sbm{A}t}$ is approximated by a rational matrix function $\hat{r}(\bm{A}t)$, where the rational function $\hat{r}(z)$ is chosen as the best rational approximation of the exponential function on the negative real axis $\R_-$. Let $\pi_{k,k}$ denote the set of rational functions $r_{k,k}(x) = p_k(x)/q_k(x)$ where $p_k$ and $q_k$ are polynomials of order $k$.  The {\sc cram} approximation of order $k$ is defined as the unique rational function $\hat{r}_{k,k} = \hat{p}_k(x)/\hat{q}_k(x)$ satisfying
    \beq
    \sup_{x \in \R_-}| \hat{r}_{k,k} (x) - e^{x} | = \inf_{r_{k,k} \in \pi_{k,k}}\left\{\sup_{x \in \R_-} | r_{k,k}(x) - e^{x}|\right\}\ .
    \eeq
The asymptotic convergence of this approximation on the negative real axis is remarkably fast with the convergence rate $ O(H^{-k}) $, where $H = 9.289\,025\,49\ldots $ is called the Halphen constant~\cite{gonchar}. It was recently discovered by Stahl and Schmelzer~\cite{stahl2009} that this convergence extends to compact subsets on the complex plane and also to Hankel contours in $\C\,\backslash\,\R_{-}$. The application of this approximation to computing the matrix exponential was originally made famous by Cody, Meinardus, and Varga in 1969 in the context of rational approximation of $e^{-x}$ on $[0, \infty)$ and it was recently resurfaced by Trefethen, Weideman, and Schmelzer~\cite{talbot}. For self-adjoint and negative semi-definite matrices, the method is guaranteed to yield an error bound in 2-norm that corresponds to the maximum error of the rational approximation on the negative real axis.
This has also been the main context for scientific applications~\cite{cody, etna, expokit}. Recently, the method {has} also {been} successfully applied to non-self-adjoint matrices with eigenvalues near the negative real axis~\cite{pusa2010, pusa2011}. These specific matrices arise from a reactor physics application, where the changes in nuclide concentrations due to radioactive decay and neutron-induced reactions are governed by a linear system $\bm{x}' = \bm{Ax}$ known as {the} burnup equations.

\section{Partial fraction decomposition form}

The main difficulty in using {\sc cram} for computing the matrix exponential is determining the coefficients of the rational function for a given $k$. In principle, the polynomial coefficients of $\hat{p}_k$ and $\hat{q}_k$ can be computed with Remez-type methods but this requires delicate algorithms combined with high-precision arithmetics. Fortunately, these coefficients have been computed to a high accuracy by Carpenter et al.\ for approximation orders $k= 0, 1, \ldots, 30$ and they are provided in~\cite{carpenter}. In practical applications, however, it is usually advantageous to employ {\sc cram} in the partial fraction decomposition ({\sc pfd}) form. For simple poles, this composition takes the form
    \beq
    \hat{r}_{k,k}(z) = \alpha_0 + \sum_{j=1}^k \frac{\alpha_j}{z-\theta_j}\ ,
    \eeq
where $ \alpha_0$ is the limit of the function $\hat{r}_{k,k}$ at infinity, and $\alpha_j$ are the residues at the poles $\theta_j$:
    \beq
    \alpha_j = \frac{\hat{p}_k(\theta_j)}{\hat{q}_k'(\theta_j)}\ .
    \label{eq: residues}
    \eeq
Since the coefficients of $\hat{r}_{k,k}$ are real, its poles form conjugate pairs, so the computational cost can be reduced to half for a real variable $x$
    \beq
    \hat{r}_{k,k}(x) = \alpha_0 + 2\, \Re\left( \sum_{j=1}^{k/2} \frac{\alpha_j}{x - \theta_j}\right)
    \label{eq: rat fct pfd real x }
    \eeq
and the matrix exponential solution may be approximated as
    \beq
    \bm{x} = \alpha_0 \bm{x}_0 + 2\, \Re\left( \sum_{j=1}^{k/2}  {\alpha_j}(\bm{A}t - \theta_j\bm{I})^{-1}\bm{x}_0 \right)\
    \label{eq: rat apprx solution for system}
    \eeq
for a real matrix $\bm{A} \in \R^{n \times n}$.

    \bfig
    \btab{c}
    \epsfig{file=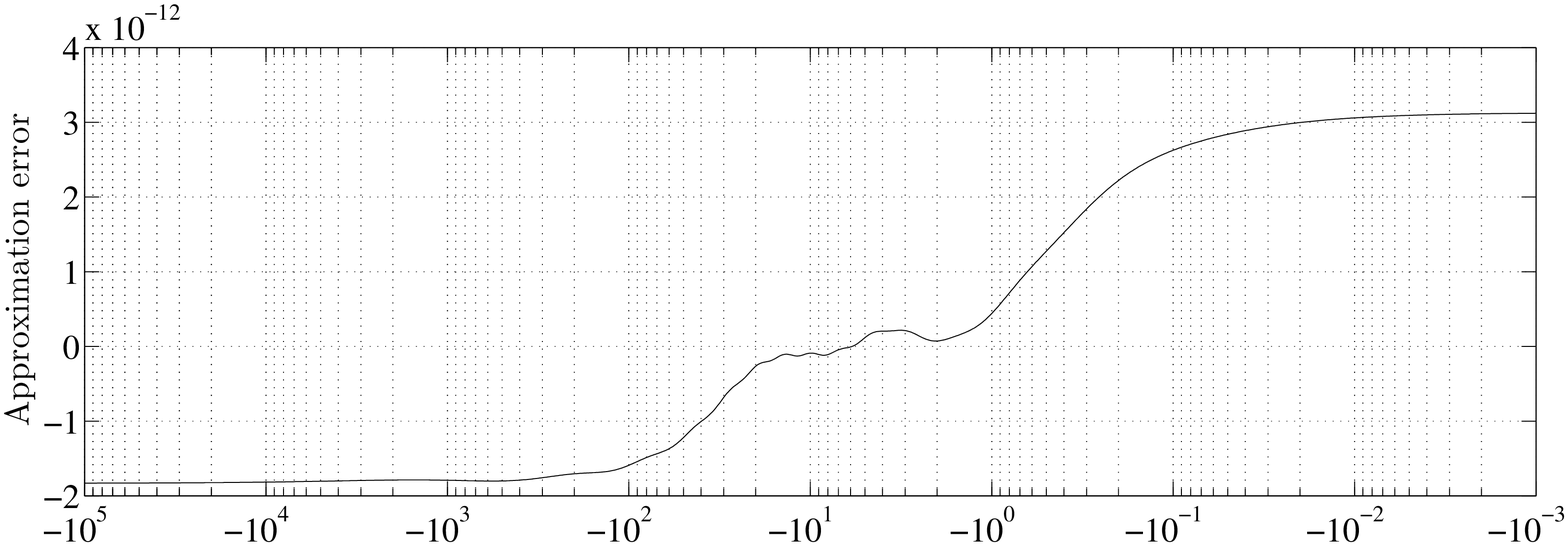,width=12cm,clip=}\\
    (a)\\[2ex]
    \epsfig{file=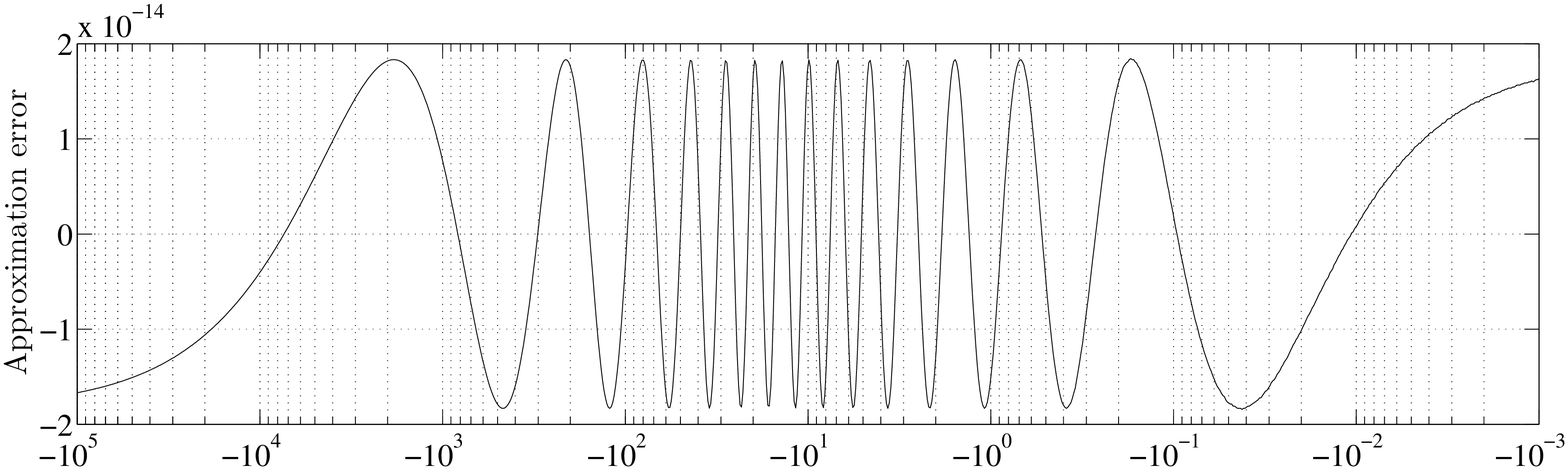,width=12cm,clip=}\\
    (b)
    \etab
    \caption{\small (a) Plot of \,$e^x - \check{r}_{14,14}(x)$\, on the negative real axis with $\check{r}_{14,14}$ computed based on the partial fraction coefficients from~\cite{gallopoulos}, (b) Plot of \,$e^x - \hat{r}_{14,14}(x)$\, based on the polynomial coefficients from~\cite{carpenter}. The plots were computed using high-precision arithmetics with 32 digits.}
    \label{fig: cram_neg_real_axis}
    \efig

Although the {\sc pfd} coefficients can in principle be computed from the polynomial coefficients, the computation of the polynomial roots may be ill-conditioned and requires great care.
The {\sc pfd} coefficients for approximation orders 10 and 14 have been provided in~\cite{gallopoulos}, and the {given} coefficients for $k=14$ have been used in several applications including the matrix exponential computing package {\sc expokit}~\cite{expokit} and {the} reactor physics code Serpent~\cite{Serpent}. However, in the latter context, it was recently observed that these reported coefficients contain errors and do not correspond to the true best approximation~\cite{pusa2011}. To illustrate this, Figure~\ref{fig: cram_neg_real_axis} shows the error of order 14 approximation on the negative real axis computed using two different sets of coefficients: the partial fraction coefficients from~\cite{gallopoulos}, {with the corresponding approximation denoted by $\check{r}_{14,14}$}, and the polynomial coefficients from~\cite{carpenter}, {with the corresponding approximation denoted by $\hat{r}_{14,14}$}. According to theory, a necessary and sufficient condition for the best approximation is that the corresponding error function equioscillates, \mbox{i.e.\ there} exists a set of points where it attains its maximum absolute value with alternating signs. Notice that the approximation computed with the coefficients from~\cite{gallopoulos} does not exhibit this behavior and in addition result{s} in a $10^2$ times poorer accuracy than expected by theory.

After discovering the erroneous behavior induced by the coefficients from~\cite{gallopoulos}, partial fraction coefficients for approximation orders $k=14$ and $k=16$ were computed from the polynomial coefficients provided in~\cite{carpenter} and {subsequently} reported in~\cite{pusa2011}. The computed {\sc pfd} coefficients are repeated here in Tables~\ref{tab: CRAM PFD 14} and~\ref{tab: CRAM PFD 16}. The computations were performed with {\sc Matlab}'s Symbolic Toolbox using high precision arithmetics with 200 digits to
ensure a sufficient accuracy.  {In Tables~\ref{tab: CRAM PFD 14} and~\ref{tab: CRAM PFD 16} the coefficients have been rounded off to 20 digits. The coefficients in~\cite{carpenter} have been also given with 20 digits' accuracy, and based on our experience, the approximation order $k=16$ is the highest for which this accuracy is sufficient for computing the {\sc pfd} coefficients. For lower approximation orders, $1 \leq k \leq 13$, the {\sc pfd} coefficients can be accurately computed with the
approximative Carath\'{e}odory--Fej\'{e}r method and a {\sc Matlab} script is provided for this purpose in~\cite{schmeltzer}.}

    \setlength\tabcolsep{5pt}
    \bt
    \caption{\small Recomputed values for the partial fraction decomposition coefficients for {\scriptsize CRAM} approximation of order 14.}
    \label{tab: CRAM PFD 14}
    \footnotesize
    \btab{c|l|l}
    \toprule
    Coefficient & \multicolumn{1}{c|}{Real part} & \multicolumn{1}{c}{Imaginary part}\\
    \midrule
    $\theta_1$ & $-8.897\,773\,186\,468\,888\,819\,9 \times 10^{0}$ & $+1.663\,098\,261\,990\,208\,530\,4 \times 10^{1}$\\[0.33ex]
    $\theta_2$ & $-3.703\,275\,049\,423\,448\,060\,3 \times 10^{0}$ & $+1.365\,637\,187\,148\,326\,817\,1 \times 10^{1}$\\[0.33ex]
    $\theta_3$ & $-0.208\,758\,638\,250\,130\,125\,1 \times 10^{0}$ & $+1.099\,126\,056\,190\,126\,091\,3 \times 10^{1}$\\[0.33ex]
    $\theta_4$ & $+3.993\,369\,710\,578\,568\,519\,4 \times 10^{0}$ & $+6.004\,831\,642\,235\,037\,317\,8 \times 10^{0}$\\[0.33ex]
    $\theta_5$ & $+5.089\,345\,060\,580\,624\,506\,6 \times 10^{0}$ & $+3.588\,824\,029\,027\,006\,510\,2 \times 10^{0}$\\[0.33ex]
    $\theta_6$ & $+5.623\,142\,572\,745\,977\,124\,8 \times 10^{0}$ & $+1.194\,069\,046\,343\,966\,976\,6 \times 10^{0}$\\[0.33ex]
    $\theta_7$ & $+2.269\,783\,829\,231\,112\,709\,7 \times 10^{0}$ & $+8.461\,737\,973\,040\,221\,401\,9 \times 10^{0}$\\[0.33ex]
    \midrule
    $\alpha_1$ & $-7.154\,288\,063\,589\,067\,285\,3 \times 10^{-5}$ & $+1.436\,104\,334\,954\,130\,011\,1 \times 10^{-4}$ \\[0.33ex]
    $\alpha_2$ & $+9.439\,025\,310\,736\,168\,877\,9 \times 10^{-3}$ & $-1.718\,479\,195\,848\,301\,751\,1 \times 10^{-2}$ \\[0.33ex]
    $\alpha_3$ & $-3.763\,600\,387\,822\,696\,871\,7 \times 10^{-1}$ & $+3.351\,834\,702\,945\,010\,421\,4 \times 10^{-1}$ \\[0.33ex]
    $\alpha_4$ & $-2.349\,823\,209\,108\,270\,119\,1 \times 10^{1}$ &  $-5.808\,359\,129\,714\,207\,400\,4 \times 10^{0}$ \\[0.33ex]
    $\alpha_5$ & $+4.693\,327\,448\,883\,129\,304\,7 \times 10^{1}$ &  $+4.564\,364\,976\,882\,776\,079\,1 \times 10^{1}$ \\[0.33ex]
    $\alpha_6$ & $-2.787\,516\,194\,014\,564\,646\,8 \times 10^{1}$ &  $-1.021\,473\,399\,905\,645\,143\,4 \times 10^{2}$ \\[0.33ex]
    $\alpha_7$ & $+4.807\,112\,098\,832\,508\,890\,7 \times 10^{0}$ &  $-1.320\,979\,383\,742\,872\,388\,1 \times 10^{0}$ \\[0.33ex]
    \midrule
    $\alpha_0$ & $+1.832\,174\,378\,254\,041\,275\,1 \times 10^{-14}$ & $+0.000\,000\,000\,000\,000\,000\,0 \times 10^{0}$\\
    \bottomrule
    \etab
    \et

    \setlength\tabcolsep{5pt}
    \bt
    \vspace*{4em}
    \caption{\small Computed values for the partial fraction decomposition coefficients for {\scriptsize CRAM} approximation of order 16.}
    \label{tab: CRAM PFD 16}
    \footnotesize
    \btab{c|l|l}
    \toprule
    Coefficient & \multicolumn{1}{c|}{Real part} & \multicolumn{1}{c}{Imaginary part}\\
    \midrule
    $\theta_1$ & $-1.084\,391\,707\,869\,698\,802\,6 \times 10^{1}$ & $+1.927\,744\,616\,718\,165\,228\,4 \times 10^{1}$\\[0.33ex] 	
    $\theta_2$ & $-5.264\,971\,343\,442\,646\,889\,5 \times 10^{0}$ & $+1.622\,022\,147\,316\,792\,730\,5 \times 10^{1}$\\[0.33ex] 	
    $\theta_3$ & $+5.948\,152\,268\,951\,177\,480\,8 \times 10^{0}$ & $+3.587\,457\,362\,018\,322\,282\,9 \times 10^{0}$\\[0.33ex]
    $\theta_4$ & $+3.509\,103\,608\,414\,918\,097\,4 \times 10^{0}$ & $+8.436\,198\,985\,884\,375\,082\,6 \times 10^{0}$\\[0.33ex]
    $\theta_5$ & $+6.416\,177\,699\,099\,434\,192\,3 \times 10^{0}$ & $+1.194\,122\,393\,370\,138\,687\,4 \times 10^{0}$\\[0.33ex]
    $\theta_6$ & $+1.419\,375\,897\,185\,665\,978\,6 \times 10^{0}$ & $+1.092\,536\,348\,449\,672\,258\,5 \times 10^{1}$\\[0.33ex]
    $\theta_7$ & $+4.993\,174\,737\,717\,996\,399\,1 \times 10^{0}$ & $+5.996\,881\,713\,603\,942\,226\,0 \times 10^{0}$\\[0.33ex]
    $\theta_8$ & $-1.413\,928\,462\,488\,886\,211\,4 \times 10^{0}$ & $+1.349\,772\,569\,889\,274\,538\,9 \times 10^{1}$\\[0.33ex]
    \midrule
    $\alpha_1$ & $-5.090\,152\,186\,522\,491\,565\,0 \times 10^{-7}$ & $-2.422\,001\,765\,285\,228\,797\,0 \times 10^{-5}$\\[0.33ex]   	
    $\alpha_2$ & $+2.115\,174\,218\,246\,603\,090\,7 \times 10^{-4}$ & $+4.389\,296\,964\,738\,067\,391\,8 \times 10^{-3}$\\[0.33ex]
    $\alpha_3$ & $+1.133\,977\,517\,848\,393\,052\,7 \times 10^{2}$ &  $+1.019\,472\,170\,421\,585\,645\,0 \times 10^{2}$\\[0.33ex]
    $\alpha_4$ & $+1.505\,958\,527\,002\,346\,752\,8 \times 10^{1}$ &  $-5.751\,405\,277\,642\,181\,997\,9 \times 10^{0}$\\[0.33ex]
    $\alpha_5$ & $-6.450\,087\,802\,553\,964\,659\,5 \times 10^{1}$ &  $-2.245\,944\,076\,265\,209\,605\,6 \times 10^{2}$\\[0.33ex]
    $\alpha_6$ & $-1.479\,300\,711\,355\,799\,971\,8 \times 10^{0}$ &  $+1.768\,658\,832\,378\,293\,790\,6 \times 10^{0}$\\[0.33ex]
    $\alpha_7$ & $-6.251\,839\,246\,320\,791\,889\,2 \times 10^{1}$ &  $-1.119\,039\,109\,428\,322\,848\,0 \times 10^{1}$\\[0.33ex]
    $\alpha_8$ & $+4.102\,313\,683\,541\,002\,127\,3 \times 10^{-2}$ & $-1.574\,346\,617\,345\,546\,819\,1 \times 10^{-1}$\\[0.33ex]
    \midrule
    $\alpha_0$ & $+2.124\,853\,710\,495\,223\,748\,8 \times 10^{-16}$ & $+0.000\,000\,000\,000\,000\,000\,0 \times 10^{0}$\\
    \bottomrule
    \etab
    \et

\clearpage
\section{Analysis of inaccurate pfd coefficients for {\footnotesize $k = 14$}}

To analyze the effect of inaccurate {\sc pfd} coefficients denoted by $\{\widetilde{\theta_j}\}$ and $\{\widetilde{\alpha}_j\}$, let $\widetilde{r}$ denote the corresponding rational approximation. The error caused by the inaccuracies in the {\sc pfd} coefficients may be estimated
        \beq
        |\hat{r}_{k,k}(z) - \widetilde{r}_{k,k}(z)| \lesssim |\alpha_0 - \widetilde{\alpha}_0| + \sum_{j=1}^k \frac{|\alpha_j|}{|z - \theta_j|^2}\, |\theta_j - \widetilde{\theta}_j|
        + \frac{1}{|z - \theta_j|}\, |\alpha_j - \widetilde{\alpha}_j|\ ,
        \label{eq: estimated error}
        \eeq
indicating that the error is the greatest in the vicinity of the poles.  It can also be seen from Eq.~\eqref{eq: estimated error} that the inaccuracy related to the poles has a greater impact near the poles, whereas the error related to the residues should begin to dominate the total error farther away from the poles. By comparing the old and the recomputed {\sc pfd} coefficients for $k=14$, it can be seen that the poles all agree to about 6 digits whereas the residues agree to about 5 digits.\,\footnote{Notice that the {\sc pfd} coefficients in~\cite{gallopoulos} are given for the rational approximation of $e^{-x}$ on $[0, \infty)$ and that they have been multiplied by a factor of two making Eq.~(37) in~\cite{gallopoulos} equivalent to Eq.~\eqref{eq: rat apprx solution for system}.} The most dramatic discrepancy occurs for the coefficient $\alpha_0$ for which the significands agree to 5 digits but the exponent value given in~\cite{gallopoulos} is $-12$, although the correct value is $-14$.

    \bfig
    \epsfig{file=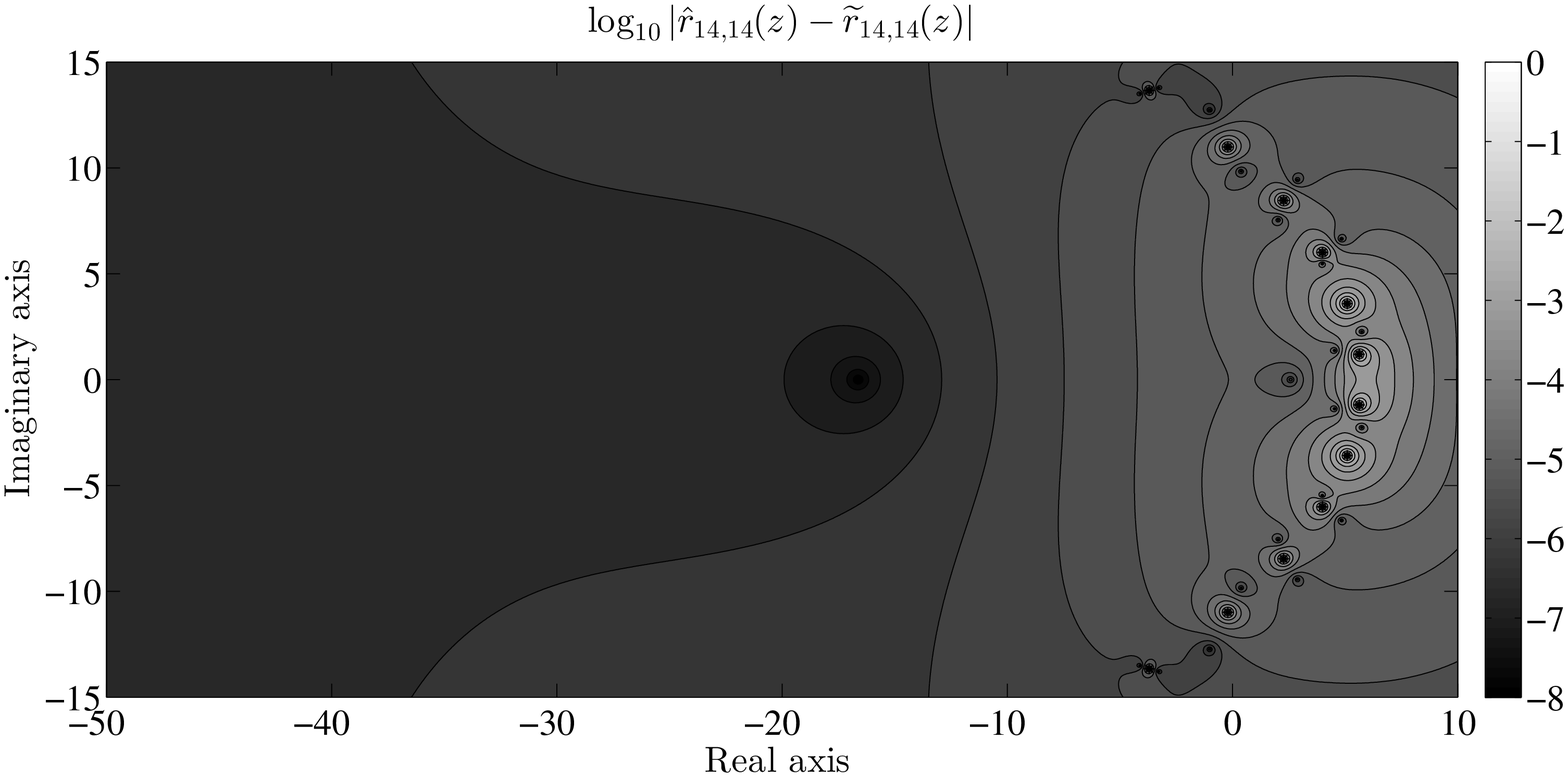,width=12cm}
    \caption{\small Plot of \,$\log_{10} |\hat{r}_{14,14}(z) - \widetilde{r}_{14,14}(z)|$\, in the complex plane. $\hat{r}_{14,14}$ was computed based on the partial fraction coefficients from Table~\ref{tab: CRAM PFD 14} and $\widetilde{r}_{14,14}$ was formed by truncating these coefficients to 6 significant digits. The poles of $\hat{r}_{14,14}$ have been marked to the plot with asterisks.}
    \label{fig: digits 6}
    \efig

On the grounds of Eq.~\eqref{eq: estimated error}, it can be estimated that coefficients with 6 correct digits should produce a rational function whose deviation from $\hat{r}_{14, 14}(z)$ is at most of the order of $10^{-3}$ on the negative real axis. Figure~\ref{fig: digits 6} shows the difference between $\hat{r}_{14, 14}(z)$ and a rational function $\widetilde{r}_{14, 14}(z)$ that was formed by truncating the {\sc pfd} coefficients of $\hat{r}_{14, 14}$ to 6 significant digits. Interestingly, as can be seen from Fig.~\ref{fig: cram_neg_real_axis}a, the approximation $\check{r}_{14,14}(x)$, corresponding to the {\sc pfd} coefficients from~\cite{gallopoulos}, yields a significantly better accuracy of order $10^{-12}$ than is expected based on the accuracy of the coefficients alone.

To investigate the matter further, let us now take the poles $\{\check{\theta}_j\}$ reported in~\cite{gallopoulos} as a starting point for constructing a rational approximation of order 14. The poles $\{\check{\theta}_j\}$ define a polynomial
    \beq
    \check{q}_{14}(x) = \prod_{j=1}^{14} (x -\check{\theta}_j )
    \eeq
whose values agree to about 6 digits with the values of the correct polynomial $\hat{q}_{14}(x)$ on the negative real axis.  The residues at the poles $\{\check{\theta}_j\}$ cannot be computed in a fully consistent manner, since the poles do not correspond to the true zeros of $\hat{q}_{14}$.  However, two alternative approaches for computing the residues can be considered. One possibility is to use the correct rational function $\hat{r}_{14, 14}$ and Eq.~\eqref{eq: residues} to compute the residues, but this is inconsistent as Eq.~\eqref{eq: residues} only holds at the true poles. Another option is to define a new rational function using $\check{q}_{14}$ as the denominator and the correct polynomial $\hat{p}_{14}$ as the numerator, after which the residues can be computed exactly using symbolic arithmetics. With both of these approaches we obtain a rational approximation, whose accuracy is of the order of $10^{-6}$ on the negative real axis. It is also worth mentioning that forming the rational function based on the poles $\{\check{\theta}_j\}$ and the correct residues $\{\alpha_j\}$ from Table~\ref{tab: CRAM PFD 14} yields an approximation whose accuracy is of the order of $10^{-7}$ on the negative real axis.

{The article~\cite{gallopoulos} by Gallopoulos and Saad does not indicate, how the reported {\sc pfd} coefficients were computed, but based on the observations regarding the accuracy of the resulting approximation}, it is evident that the values given for the residues somehow compensate for the {inaccuracies in the} poles $\{\check{\theta}_j\}$ and it seems likely that they have been optimized to minimize the deviation from $\hat{r}_{14,14}$ on the negative real axis. In fact, using the poles $\{\check{\theta}_j\}$ and standard least squares optimization in {\sc Matlab} with $10^7$ points chosen from the interval $[-10^3,-10^{-10}]$, we were able to produce residues yielding only a slightly worse accuracy of order $10^{-11}$. In any case, it should be noted that optimizing the residues properly in the Chebyshev sense would essentially form a problem of comparable difficulty as the original problem of determining $\hat{r}_{14,14}$.

\end{document}